\documentclass[a4paper,10pt]{article}
\usepackage{tikz}
\usetikzlibrary{automata,positioning}
\usepackage{multirow}
\usepackage{pgfplots}
\usetikzlibrary{positioning}
\usetikzlibrary{arrows}
\usetikzlibrary{matrix,chains,positioning,decorations.pathreplacing,arrows}
\usetikzlibrary{calc,arrows.meta,positioning}
\usepackage{subcaption}
\usepackage{float}
\usepackage{needspace}
\newtheorem{result}{\textbf{Result}}

\newtheorem{remark}{\textbf{Remark}}

\newtheorem{definition}{\textbf{Definition}}

\usepackage{cite}
\usepackage{amsmath,amssymb,amsfonts}
\usepackage{algorithmic}
\usepackage{graphicx}
\usepackage{textcomp}
\usepackage{xcolor}
\bgroup\catcode`\#=12\egroup\newcommand{\qed}{\hfill$\blacksquare$}

\begin{document}
\title{A convex optimization approach to online set-membership EIV identification of LTV systems}

\author{S. M. Fosson${}$, D. Regruto${}$, T. Abdalla${}$ and A. Salam${}$}
\maketitle

\abstract{%
This paper addresses the problem of recursive set-membership identification for linear time varying (LTV) systems when both input and output measurements are affected by bounded additive noise. First we formulate the problem of online computation of the parameter uncertainty intervals (PUIs) in terms of nonconvex polynomial optimization. Then, we propose a convex relaxation approach based on McCormick envelopes to solve the formulated problem to the global optimum by means of linear programming. The effectiveness of the proposed identification scheme is demonstrated by means of two simulation examples.}


\section{Introduction}

Recursive parameter estimation of linear systems has continuously attracted the attention of the automatic control community in the last decades. This topic is of particular interest in the context of linear time-varying (LTV) systems, where the parameter variations need to be tracked online.

A relevant number of contributions addressing the problem of recursive identification of LTV systems can be found in the context of classical system identification, where the noise affecting the measurements is assumed to be statistically described. The interested reader can find details in the survey paper \cite{lju90} and in the book \cite{Honig,niedzwieck}.

A worthwhile alternative to the stochastic noise description, inspired by the seminal work of Schweppe \cite{sch68}, is the so-called bounded-errors or set-membership (SM) characterization, where uncertainties are assumed to belong to a given set, see, e.g., the book \cite{mil96} for an introduction to the theory. In the SM framework, all parameter vectors belonging to the feasible parameter set (FPS), i.e., parameters consistent with the measurements, the error bounds and the assumed model structure, are feasible solutions to the identification problem. The objective of any SM algorithm is either to optimally select a single solution in the FPS (pointwise SM estimators) or to compute uncertainty bounds for the parameters (set-valued SM estimators). In this work we focus our attention on this second class. A number of algorithms have been proposed to address the problem of computing parameter bounds for LTV systems. The idea common to all the approaches is to recursively approximate the FPS by means of simply-shaped sets: ellipsoids are considered in \cite{nor90}, polyedrals in \cite{pie94} and zonotopes in \cite{cha06,bra06}, while orthotopic regions have been recently considered in \cite{cas17}. 

All the aforementioned papers formulate the identification problem with reference to the equation error structure. To the best of the authors'  knowledge, the first attempt to address the problem of recursive SM identification for LTV systems in the errors-in-variables (EIV) framework, i.e., when both the input and the output measurements are affected by noise, has been presented in our previous contribution \cite{cer20}. In this paper, we show that the parameter uncertainty intervals (PUIs) can be exactly computed at each time iteration by solving a set of simple linear programming problems, provided that the sign of the parameters are a-priori known. However, if such an information is not available, the problem requires the solution of a number of computationally expensive non-convex polynomial optimization problems.     

In order to overcome this limitation, in this work we propose a different convex relaxation strategy which does not require any a-priori information on the parameters sign. The proposed approach is based on the concept of McCormick envelopes originally proposed in \cite{McCormick76}.

The paper is organized as follows. Section \ref{Problem Formulation} is devoted to the problem formulation. In Section \ref{Bounding the parameters of the LTV system}, we briefly review the results of our previous contribution \cite{cer20} in order to highlight the mathematical structure of the problem. The novel convex relaxation approach, based on McCormick envelopes, is presented in Section \ref{McCormick convexication}. The effectiveness of the proposed identification scheme is shown in Section \ref{Simulated Examples} by means of two simulation examples. Concluding remarks end the paper.
\section{Problem Formulation}\label{Problem Formulation}
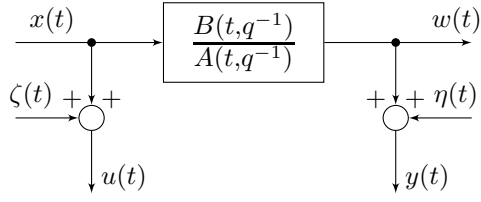
\begin{figure}[t!]
	\tikzstyle{block} = [draw,fill=white!20, rectangle,
	minimum height=3em, minimum width=6em]
	\tikzstyle{block_NL} = [draw,fill=white!20, rectangle, 
	minimum height=3em, minimum width=6em]
	\tikzstyle{sum} = [draw, fill=white!20, circle, node distance=0cm]
	\tikzstyle{CNode} = [draw, fill=black!100, circle, inner sep=0pt,minimum size=0.1cm,node distance=0cm]
	\tikzstyle{input} = [coordinate]
	\tikzstyle{output} = [coordinate]
	\tikzstyle{pinstyle} = [pin edge={to-,thick,black}]
	\begin{center}
		
		\begin{tikzpicture}[auto, node distance=1cm,>=latex']
		
		\node [input, name=input] {};
		\node [output, name=output] {};
		
		\node[CNode, right of=input,node distance=1cm, name=CNode] {};
		\node [block, right of=CNode, node distance=2cm] (system) {\Large${\frac{B(t,q^{-1})}{A(t,q^{-1})}}$};
		\node[CNode, right of=system,node distance=2cm, name=GNode] {};
		\node [sum, below of=CNode,node distance=1cm] (sum1) {};
		\node [sum, below of=GNode,node distance=1cm] (sum2) {};
		\node [output, right of=system,node distance=2.5cm] (w) {};
		\node [output, right of=GNode] (w) {};
		\node [output, right of=sum2] (eta) {};
		\node [output, left of=sum1] (zeta) {};
		\node [output, below of=sum1] (ut) {};
		\node [output, below of=sum2] (yt) {};
	 
		\draw [draw,->] (input) -- node [near start] {$x(t)$} (system);
		\draw [draw,->] (system) -- node [pos=0.9] {$w(t)$} (w);
		\draw [->] (CNode) --  node [pos=0.9] {$+$} node {} (sum1);
		\draw [->] (zeta) -- node [pos=0.9] {$+$} node [near start] {$\zeta(t)$} (sum1);
		\draw [->] (sum1) -- node [near end] {${u}(t)$} (ut);
		\draw [->] (GNode) --  node [pos=0.9,left] {$+$} node {} (sum2);
		\draw [->] (eta) -- node [pos=0.9,above] {$+$} node [near start, above] {$\eta(t)$} (sum2);
		\draw [->] (sum2) -- node [near end] {$y(t)$} (yt);
		
		\end{tikzpicture}
	\end{center}
	\caption{EIV model structure}
	\label{fig:EIV Structure}
\end{figure}

Let us consider the SISO discrete-time LTV system, depicted in Fig. \ref{fig:EIV Structure}, described in terms of the following linear difference equation
\begin{equation}\label{eq 2}
        A(t,q^{-1})w(t)=B(t,q^{-1})x(t),
\end{equation}
where $x(t)$ and $w(t)$ are the noise free input and output signals respectively, and $A(t,\cdot)$, $B(t,\cdot)$ are polynomials in the backward shift operator $q^{-1} (i.e., q^{-1}u(t) = u(t-1)$) given by
\begin{equation}
    A(t,q^{-1})=1+a_{1}(t)q^{-1}+\dots+a_{n_a}(t)q^{-n_{a}},
        \label{Polynomial A}
\end{equation}

\begin{equation}
    B(t,q^{-1})=b_{0}(t)+b_{1}(t)q^{-1}+\dots+b_{n_b}(t)q^{-n_{b}},
    \label{Polynomial B}
\end{equation}
where $n_a \geq n_b$.\\
The unknown parameter vector $\theta(t) \in {\mathbb{R}}^{n_p}$ to be estimated at each time instant $t$ is
\begin{equation}
    \theta(t)=[a_{1}(t),\dots,a_{n_a}(t),
    b_{0}(t),\dots,b_{n_b}(t)]^T,
     \label{theta vector}
     \end{equation}
where $n_p=n_a+n_b+1$.
\\At each generic time instant $t$, the value of $k$-th component $\theta_k(t)$ of the parameter vector $\theta(t)$  is described as 
\begin{equation}
    \theta_k(t)=\theta_k(t-1)+\delta_{\theta_k}(t), \quad \quad k = 1,\dots,n_p,
    \label{theta variation}
\end{equation}
where $\delta_{\theta_k}(t)$ is the parameter variation between two consecutive generic time instants $t-1$ and $t$, assumed to be unknown but bounded, i.e.,
\begin{equation}
    |\delta_{\theta_k}(t)| \leq \Delta_{\theta_k}, \forall t
    \end{equation}
where $\Delta_{\theta_k}$, $k=1,\ldots,n_p$, are known variation bounds.\\
Both input and output data are corrupted by additive noise $\zeta(t)$ and $\eta(t)$ respectively,
\begin{equation}\label{eq: input}
    u(t) = x(t)+\zeta(t),
\end{equation}

\begin{equation}\label{eq: output}
    y(t) = w(t)+\eta(t).
\end{equation}
Each sample of the noise sequences $\zeta(t)$ and $\eta(t)$ is bounded by known constants $\Delta_{\zeta}$ and $\Delta_{\eta}$, i.e.
\begin{equation}
    |\zeta(t)| \leq \Delta_{\zeta},  \forall t
    \label{input error bound}
\end{equation}
\begin{equation}
    |\eta(t)| \leq \Delta_{\eta},  \forall t.
    \label{output error bound}
\end{equation}
\\
According to the problem formulation presented in \cite{cer20}, the FPS, at a generic time instant $t$, is defined as,
\begin{equation}
   \begin{split}
    {\mathcal{D}}_{\theta}(t) = \{ & \theta(t) \in {\mathbb{R}}^{n_p}:A(q^{-1},t)(y(t)-\eta(t))
    =B(q^{-1},t)(u(t)-\zeta(t)),\\
     &\theta_k(t)=\theta_k(t-1)+\delta_{\theta_k}(t), \ |\delta_{\theta_k}(t)| \leq \Delta_{\theta_k},\\ &\underline{\theta}_{k}(t-1)\leq\theta_k(t-1)\leq\overline{\theta}_{k}(t-1),\\& k = 1,...,n_p, \ |\eta(t)|\leq\Delta_{\eta},    \ |\zeta(t)|\leq\Delta_{\zeta}\},
    \end{split}
    \label{FPS}
\end{equation}
where $\underline{\theta}_{k}(t-1)$ and $\overline{\theta}_{k}(t-1)$ are bounds on $\theta_k(t-1)$ computed at time $t-1$. 

In this work we address the problem of online computation of the parameter uncertainty intervals (PUIs) defined as
\begin{equation}
    PUI_{k}(t) = [\underline{\theta}_{k}(t) \ , \ \overline{\theta}_{k}(t)],
    \label{PUI}
\end{equation}
where
\begin{equation}\label{PUI1}
\begin{aligned}
    \underline{\theta}_{k}(t) = \underset{\theta(t) \in {\mathcal{D}}_{\theta}(t)}{\text{min}} \theta_{k}(t),
    \end{aligned}
\end{equation}

\begin{equation}\label{PUI2}
\begin{aligned}
    \overline{\theta}_{k}(t) = \underset{\theta(t) \in {\mathcal{D}}_{\theta}(t)}{\text{max}} \theta_{k}(t). 
    \end{aligned}
\end{equation}
 Initial bounds $\underline{\theta}_{k}(0)$ and $\overline{\theta}_{k}(0)$ are assumed to be a-priori known.
\section{Bounding the parameters of the LTV system}\label{Bounding the parameters of the LTV system}
In this section we briefly review the approach proposed in our previous contribution \cite{cer20}  for the computation of the solution to problems \eqref{PUI1} and \eqref{PUI2}.

\noindent The following result provides insight into the mathematical structure of the optimization problem to be solved. \\
\begin{result}\normalfont{\textbf{Computation of $\underline{\theta}_{j}(t)$ and $\overline{\theta}_{j}(t)$ via polynomial optimization}} 
\\
Lower and upper bounds on the $k-$th component $\theta_k(t)$ of the parameter vector can be computed solving the following (nonconvex) polynomial optimization problem :
\begin{equation}
\begin{aligned}
    & \underset{\theta(t)}{\text{min}} \ J_{\theta_k}(t)\\
    & \ s.t:
    \\
    & \left\{\begin{matrix}
                A(t,q^{-1})(y(t)-\eta(t)) = B(t,q^{-1})(u(t)-\zeta(t)), \vspace{0.2cm}\\ 
        \hspace{-4.5cm}-\Delta_\eta \leq \eta_t \leq \Delta_{\eta_t},\vspace{0.2cm}\\
         \hspace{-4.5cm}-\Delta_\zeta \leq \zeta_t \leq \Delta_{\zeta_t},\vspace{0.2cm}\\ \hspace{+0.1cm}\underline{\theta}_{k}(t-1)-\Delta_{\theta_k}(t)\leq\theta_k(t)\leq\overline{\theta}_{k}(t-1)+\Delta_{\theta_k}(t),\vspace{0.2cm}\\
         \hspace{-4.8cm}k = 1,\ldots,n_p
\end{matrix}\right.
    \end{aligned}
    \label{second OP}
\end{equation}

where $\ J_{\theta_k}(t)= {\theta_k}(t)$ for computing the lower bound $\underline{\theta}_{k}(t)$, or $\ J_{\theta_k}(t)= -{\theta_k}(t)$ for the upper bound $\overline{\theta}_{k}(t)$.\qed
\\
\end{result}

Problem \eqref{second OP} is directly derived from \eqref{PUI1} and \eqref{PUI2} by rewriting the constraints describing the FPS (11) in compact form. 
Non-convexity of optimization problem \eqref{second OP} is due to the presence of bilinear terms (involving unknown variables $\theta$, $\eta$ and $\zeta$) in the equality $A(t,q^{-1})(y(t)-\eta(t)) = B(t,q^{-1})(u(t)-\zeta(t))$. As discussed in our previous contribution \cite{cer20}, problem \eqref{second OP} can be solved to global optimum by means of linear programming if a-priori information on the sign of the parameters are available. 
In fact, under such an assumption,
the bilinear model equation can be rewritten as:
\begin{equation}
    (\varphi(t)-\Delta_{\varphi}(t))\theta\leq y_t+\Delta_\eta,
\end{equation}

\begin{equation}
    (\varphi(t)+\Delta_{\varphi}(t))\theta\geq y_t-\Delta_\eta,
\end{equation}
where $\varphi(t)$ is the measurement regressor defined as
\begin{equation}
     \varphi(t) = [-y(t-1),\ldots,-y(t-n_a),u(t),\ldots u(t-n_b)],
\end{equation}
while $\Delta_{\varphi}(t)$ is given by:
\begin{equation}
\begin{aligned}
      \Delta_{\varphi}(t) = [&\Delta_\eta sgn(a_1(t)),\ldots,\Delta_\eta sgn(a_{n_a}(t))\\&\Delta_\zeta sgn(b_0(t)),\ldots,\Delta_\zeta sgn(b_{n_b}(t))].
      \end{aligned}
\end{equation}
\begin{remark}\normalfont
It is worth noting that, in case no information is available about the sign of the parameters, the approach proposed in paper \cite{cer20} cannot be applied to convert \eqref{second OP} to a linear program. To address this drawback, on the one hand, one can resort to prior estimation of the signs, see, e.g., \cite{fox20}; nevertheless, this approach may be computationally unfeasible in online identifcation. On the other hand,  convex relaxation techniques guaranteed to converge to the global optimum of nonconvex polynomial optimization problems  are available in the literature, see, e.g., \cite{Lasserre,Parrilo,Chesi}. Nevertheless, such methods require the solution of  large-dimensional semidefinite programming problems even when the number of parameters is relatively small; therefore, they cannot be applied in the framework of online estimation of LTV systems, due to their high computational complexity in terms of both computational time and memory resources requirements. In this work, we propose a novel approach which does not require any information about the signs of the parameters.
\end{remark}
\section{McCormick envelopes based convex relaxation}\label{McCormick convexication}
In this section, we propose an approach to reformulate problem \eqref{second OP} in terms of convex optimization. The main idea is to exploit the concept of McCormick envelopes \cite{McCormick76} to replace the bilinear terms in \eqref{second OP} with a set of linear inequalities, without introducing any conservativeness.
\\
Let us first rewrite model equations \eqref{eq 2}, \eqref{eq: input} and \eqref{eq: output} in the following compact form:
\begin{equation}\label{full eq}
\begin{aligned}
        y(t)-\eta(t) = &-\sum_{i=1}^{n_a}a_i(t)y(t-i)+\sum_{j=0}^{n_b}b_j(t)u(t-j)+\\&+\sum_{i=1}^{n_a}a_i(t)\eta(t-i)-\sum_{j=0}^{n_b}b_j(t)\zeta(t-j).
        \end{aligned}
\end{equation}\\
Then, in order to eliminate the bilinear (nonconvex) terms in \eqref{full eq}, we define the following new variables
\begin{equation}\label{Mc1}
    {\mathcal{M}}_{a_i}(t) = a_i(t)\eta(t-i), \quad \quad i = 1,\ldots,n_a,
\end{equation}

\begin{equation}\label{Mc2}
    {\mathcal{M}}_{b_j}(t) = b_j(t)\zeta(t-j), \quad \quad j = 0,\ldots,n_b
\end{equation}
which allow us to rewrite \eqref{full eq} as follows
\begin{equation}\label{eq with M}
\begin{aligned}
        &y(t)+\sum_{i=1}^{n_a}a_i(t)y(t-i)-\sum_{j=0}^{n_b}b_j(t)u(t-j)\\&-\sum_{i=1}^{n_a}{\mathcal{M}}_{a_i}(t)+\sum_{j=0}^{n_b}{\mathcal{M}}_{b_j}(t) = \eta(t).
        \end{aligned}
\end{equation}
Since $\eta$ is known to be bounded according to \eqref{output error bound}, the following inequality is finally obtained:
\begin{equation} \label{eq with M bound}
\begin{aligned}
        &|y(t)+\sum_{i=1}^{n_a}a_i(t)y(t-i)-\sum_{j=0}^{n_b}b_j(t)u(t-j)-\\&-\sum_{i=1}^{n_a}{\mathcal{M}}_{a_i}(t)+\sum_{j=0}^{n_b}{\mathcal{M}}_{b_j}(t)|\leq \Delta_\eta.
        \end{aligned}
\end{equation}
\\
Upper and lower bounds on ${\mathcal{M}}_{a_i}(t)$ and ${\mathcal{M}}_{b_j}(t)$, can be obtained by relying on the concept of McCormick envelopes. \\

{\begin{definition}\normalfont{(\textbf{McCormick envelopes \cite{McCormick76}})}\label{def mc}\\
Given two bounded variables $x, y \in {\mathbb{R}}$, $x^{LB}\leq x \leq x^{UB}$ and $y^{LB}\leq y \leq y^{UB}$, the product $w = xy$, satisfies the following inequalities:
\begin{equation}
    w \geq x^{LB}y+xy^{LB}-x^{LB}y^{LB},
\end{equation}

\begin{equation}
    w \geq x^{UB}y+xy^{UB}-x^{UB}y^{UB},
\end{equation}

\begin{equation}
    w \leq x^{UB}y+xy^{LB}-x^{UB}y^{LB},
\end{equation}

\begin{equation}
    w \leq xy^{UB}+x^{LB}y-x^{LB}y^{UB}.
\end{equation}\qed
\end{definition}} 
Direct application of the concept of McCormick envelopes in Definition 1 to equations \eqref{Mc1} and \eqref{Mc2} leads to the following result.
\begin{result}\normalfont{\textbf{Computation of  ${\mathcal{M}}_{\theta_k}(t)$ bounds}}\label{M bound}\\
Let ${\mathcal{M}}_{\theta_k}(t)$ be the generic term in either equation \eqref{Mc1} or equation \eqref{Mc2}. ${\mathcal{M}}_{\theta_k}(t)$ satisfies the following set of inequalities
\begin{equation}\label{first}
    {\mathcal{M}}_{\theta_k}(t) \geq \theta_k^{LB}(t)\epsilon(t-\lambda)-\theta_k(t)\Delta_\epsilon+\theta_k^{LB}(t)\Delta_\epsilon,
\end{equation}

\begin{equation}
    {\mathcal{M}}_{\theta_k}(t) \geq \theta_k^{UB}(t)\epsilon(t-\lambda)+\theta_k(t)\Delta_\epsilon-\theta_k^{UB}(t)\Delta_\epsilon,
\end{equation}

\begin{equation}
    {\mathcal{M}}_{\theta_k}(t) \leq \theta_k^{UB}(t)\epsilon(t-\lambda)-\theta_k(t)\Delta_\epsilon+\theta_k^{UB}(t)\Delta_\epsilon,
\end{equation}

\begin{equation}\label{last}
    {\mathcal{M}}_{\theta_k}(t) \leq \theta_k(t)\Delta_\epsilon+\theta_k^{LB}(t)\epsilon(t-\lambda)-\theta_k^{LB}(t)\Delta_\epsilon,
\end{equation}

where $\epsilon(t)=\eta(t)$ and  $\Delta_\epsilon=\Delta_\eta$ for  ${\mathcal{M}}_{\theta_k}(t)={\mathcal{M}}_{a_i}(t)$ in equation \eqref{Mc1}, while $\epsilon(t)=\zeta(t)$ and $\Delta_\epsilon=\Delta_\zeta$ for ${\mathcal{M}}_{\theta_k}(t)={\mathcal{M}}_{b_j}(t)$ in equation \eqref{Mc2}. Bounds $\theta_k^{LB}(t)$ and $\theta_k^{UB}(t)$ are given by
\begin{equation}
    \theta_k^{LB}(t) = \underline{\theta}_k(t-1)-\Delta_{\theta_k},
\end{equation}

\begin{equation}
    \theta_k^{UB}(t) = \overline{\theta}_k(t-1)+\Delta_{\theta_k}.
\end{equation}
\\
Variable $\lambda=i$ for $\theta_k=a_i$ in equation \eqref{Mc1}, while $\lambda=j$ for $\theta_k=b_j$ in equation \eqref{Mc2}.
\end{result}\qed
\\
Thanks to Result \ref{M bound}, we are now in the position to state the main result of the paper.
\\
\begin{result}\normalfont{\textbf{Computation of PUIs by means of linear programming}\label{PUI_LP}
\\
The global optimal solution to optimization problem \eqref{second OP} can be computed solving the following linear program:}
\begin{equation}\label{final OP}
\begin{aligned}
    & \underset{\theta(t)}{\text{min}} \ J_{\theta_k}(t)\\
    & \ s.t:
    \\
    & \left\{\begin{matrix}
    \begin{aligned}
            &|y(t)+\sum_{i=1}^{n_a}a_i(t)y(t-i)-\sum_{j=0}^{n_b}b_j(t)u(t-j)\\&-\sum_{i=1}^{n_a}{\mathcal{M}}_{a_i}(t)+\sum_{j=0}^{n_b}{\mathcal{M}}_{b_j}(t)|\leq \Delta_\eta, \vspace{0.2cm}\\ 
        &{\mathcal{M}}_{\theta_k}(t) \geq \theta_k^{LB}(t)\epsilon(t-\lambda)-\theta_k(t)\Delta_\epsilon+\theta_k^{LB}(t)\Delta_\epsilon,\vspace{0.2cm}\\ 
        &{\mathcal{M}}_{\theta_k}(t) \geq \theta_k^{UB}(t)\epsilon(t-\lambda)+\theta_k(t)\Delta_\epsilon-\theta_k^{UB}(t)\Delta_\epsilon,\vspace{0.2cm}\\ 
        &{\mathcal{M}}_{\theta_k}(t) \leq \theta_k^{UB}(t)\epsilon(t-\lambda)-\theta_k(t)\Delta_\epsilon+\theta_k^{UB}(t)\Delta_\epsilon,\vspace{0.2cm}\\ 
        &{\mathcal{M}}_{\theta_k}(t) \leq \theta_k(t)\Delta_\epsilon+\theta_k^{LB}(t)\epsilon(t-\lambda)-\theta_k^{LB}(t)\Delta_\epsilon,\vspace{0.2cm}\\
        &\theta_k^{LB}(t) = \underline{\theta}_k(t-1)-\Delta_{\theta_k},\vspace{0.2cm}\\
        &\theta_k^{UB}(t) = \overline{\theta}_k(t-1)+\Delta_{\theta_k},\vspace{0.2cm}\\
        &\theta_k^{LB}(t)\leq\theta_k(t)\leq\theta_k^{UB}(t),\vspace{0.2cm}\vspace{0.2cm}\\
        &k = 1,\ldots,n_p
             \end{aligned}
\end{matrix}\right.
    \end{aligned}
\end{equation}
\qed
\end{result}
Result \ref{PUI_LP} is proved by replacing the first constraint in \eqref{second OP} ($A(t,q^{-1})(y(t)-\eta(t)) = B(t,q^{-1})(u(t)-\zeta(t))$) with the linear inequality \eqref{eq with M bound} and the bounds on ${\mathcal{M}}_{\theta_k}$ defined in Result \ref{M bound}.

\section{Simulation examples} \label{Simulated Examples}

In order to demonstrate the effectiveness of the proposed approach, two numerical examples are presented in this section.
Computations are performed on an Intel Core i7-10510U $@$ 1.80GHz computer with 16 GB RAM, using IBM ILOG CPLEX optimizer under Matlab R2018b.

\subsection{Example 1}
\label{ex. 1}
Let us consider the first order LTV system described by the following input-output equation,
\begin{equation}
    w(t) = -a_1(t)w(t-1)+b_1(t)x(t-1),
\end{equation}
where,
\begin{align*}
     b_{1}(t) & = -2+0.5\sin{\left[\frac{2\pi t}{750}\right]},\\
     a_{1}(t) & = 0.2+0.4\sin{\left[\frac{2\pi t}{500}\right]},
\end{align*}
and parameter variation bounds are $\Delta_{b_1}=\pi/750$ and $\Delta_{a_1}=0.8\pi/500$.

\begin{figure*}[t!]
  \begin{minipage}{0.43\textwidth}
\includegraphics[width=1\textwidth]{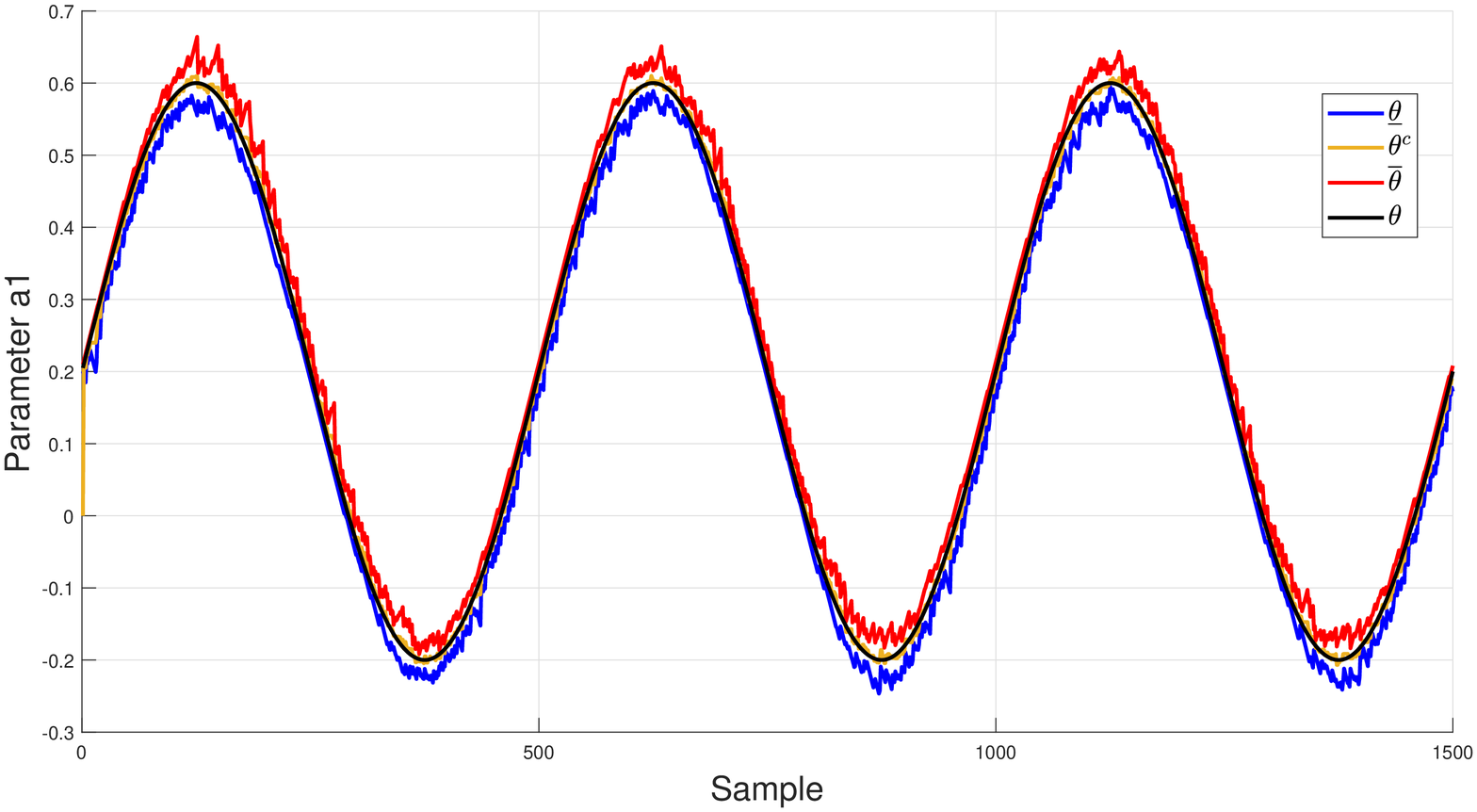}\\
          \includegraphics[width=1\textwidth]{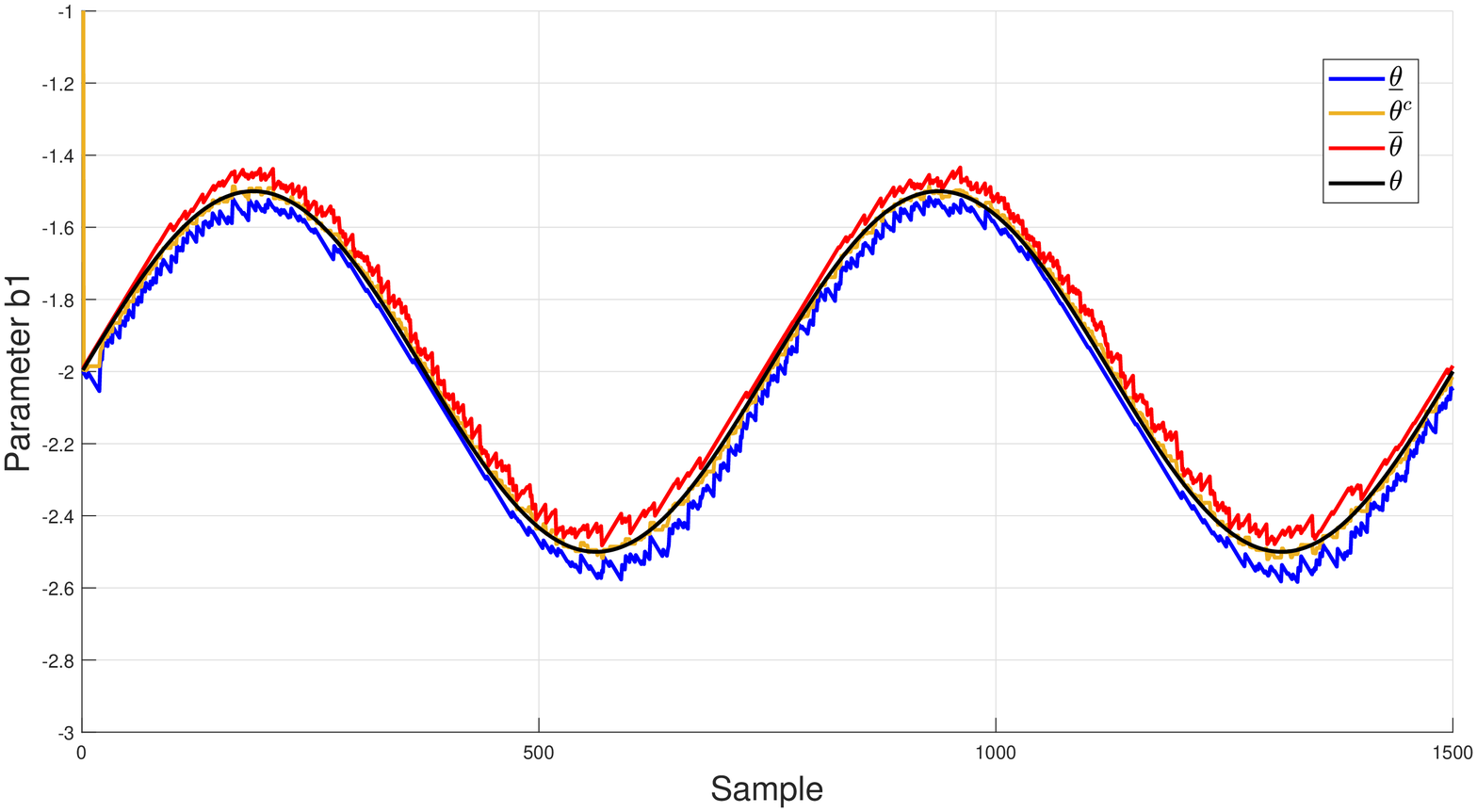}
          \subcaption{$SNR_x = 47$ dB and $SNR_w = 46$.}
          \end{minipage}
           \quad
\begin{minipage}{0.43\textwidth}
\includegraphics[width=1\textwidth]{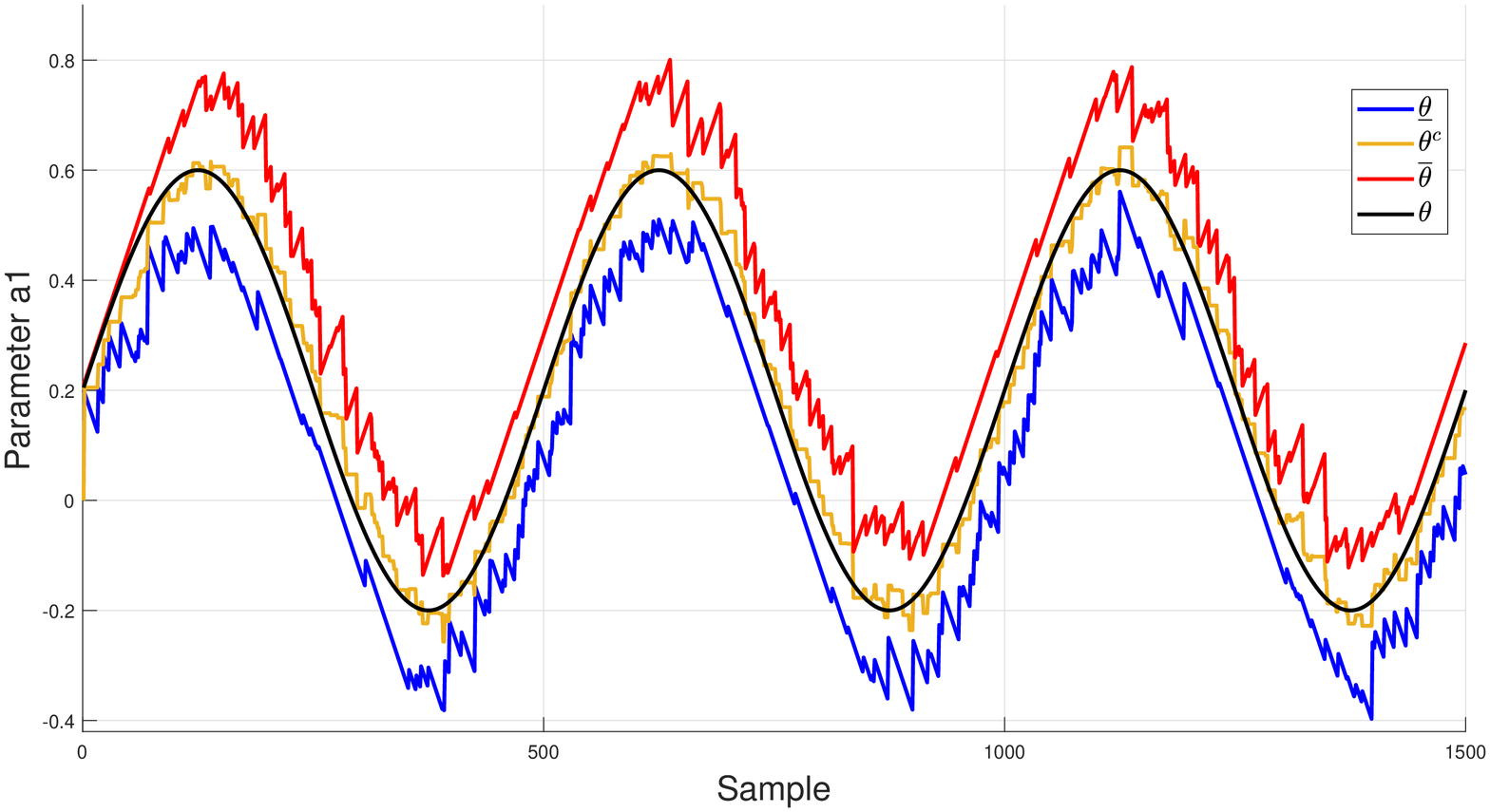}\\
          \includegraphics[width=1\textwidth]{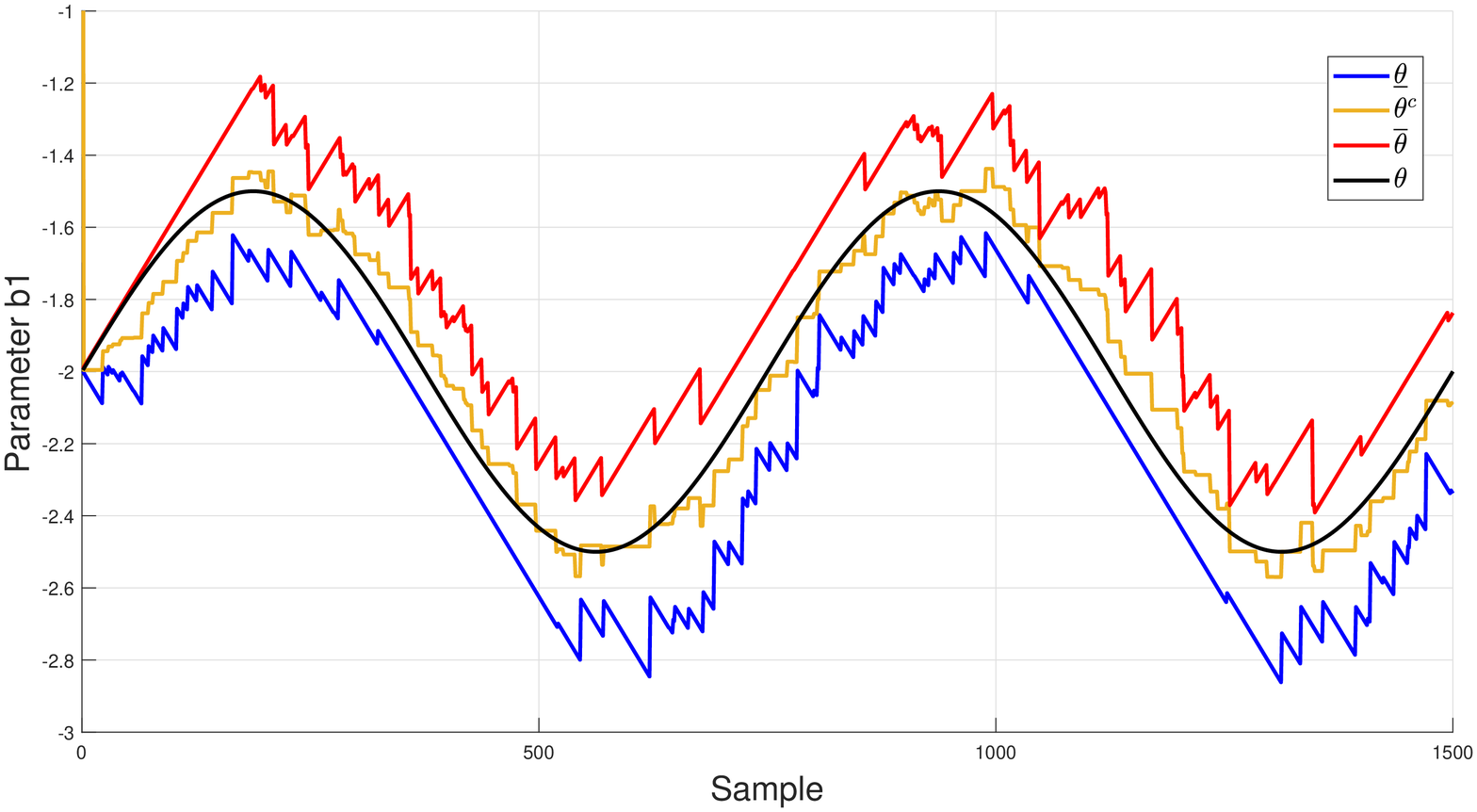}
          \subcaption{$SNR_x = 27$ dB and $SNR_w = 26$.}
          \end{minipage} 
          \caption{Example 1: Computed PUIs and central estimate through the proposed online identification scheme.}
          \label{fig: 51}
\end{figure*}

 The input is a random sequence uniformly distributed between $[-1 \ , \ +1]$. Both input and output sequence are corrupted by random additive noise, uniformly distributed between $[-\Delta_{\zeta} \ , \ \Delta_{\zeta}]$ and $[-\Delta_{\eta} \ , \ \Delta_{\eta}]$, respectively. The error bounds $\Delta_{\zeta}$ and $\Delta_{\eta}$ are chosen in such a way as to simulate two different values for both the input ($SNR_x = [47, 27]$) and the output ($SNR_w = [46, 26]$) signal-to-noise ratios, respectively defined as:
\begin{equation}
    SNR_x = 10\log\left \{\sum_{t=1}^{N}x_{t}^{2}\bigg/ \sum_{t=1}^{N}\zeta_{t}^{2}  \right \},
    \label{eq:28}
\end{equation}
\begin{equation}
    SNR_w = 10\log\left \{\sum_{t=1}^{N}w_{t}^{2}\bigg/ \sum_{t=1}^{N}\eta_{t}^{2}  \right \},
    \label{eq:29}.
\end{equation}
In this example we consider a data set of length $N=1500$. Fig. \ref{fig: 51} shows the computed bounds $\underline{\theta}$ and $\overline{\theta}$, at each sampling instant, alongside the central estimates $\theta^c$ given by
\begin{equation}
    \theta^{c}_{k}(t) = \frac{\overline{\theta}_{k}(t)+\underline{\theta}_{k}(t)}{2}, \quad \quad k = 1,\ldots,n_p,
\end{equation}
which represent the Chebyshev centers in the $\ell_\infty$-norm of ${\mathcal{D}}_{\theta}(t)$ and enjoys peculiar optimality properties (see \cite{Kacewicz} for details). Average CPU time, at each recursion, is about 1.5 ms. We can clearly observe from these figures that parameter $a_1$ changes sign several times, but this has no effect on the performance of the algorithm, and the true parameter is always included in the interval between $\underline{\theta}$ and $\overline{\theta}$.

\subsection{Example 2}
This second example is taken from \cite{cer20}. The system to be identified is a second order LTV system described by the following transfer function,
\begin{equation}
    G(q^{-1},t) = \frac{b_{1}(t)q^{-2}}{1+a_{1}(t)q^{-1}+a_{2}(t)q^{-2}},
\end{equation}
where $a_2 = 0.25$ is a fixed parameters, while $a_1$ and $b_1$ vary according to,
\begin{align*}
     b_{1}(t) & = 0.8+0.3\sin{\left[\frac{2\pi t}{2000}\right]},\\
     a_{1}(t) & = \ \  1+0.1\sin{\left[\frac{2\pi t}{1000}\right]},
\end{align*}
and parameter variation bounds are $\Delta_{b_1}=0.6\pi/2000$ and $\Delta_{a_1}=0.2\pi/1000$.

\begin{figure*}[!h]
  \begin{minipage}{0.43\textwidth}
\includegraphics[width=1\textwidth]{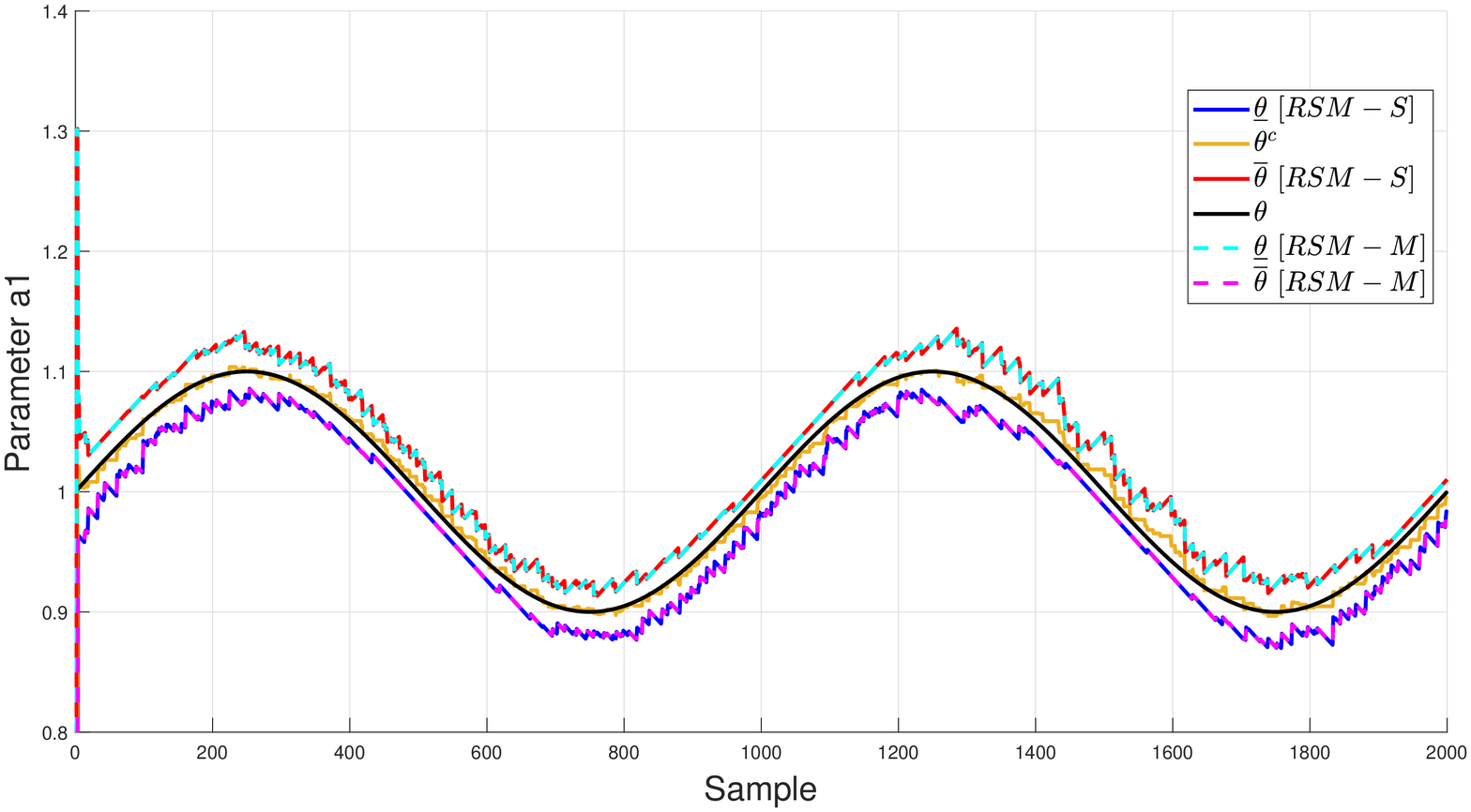}\\
          \includegraphics[width=1\textwidth]{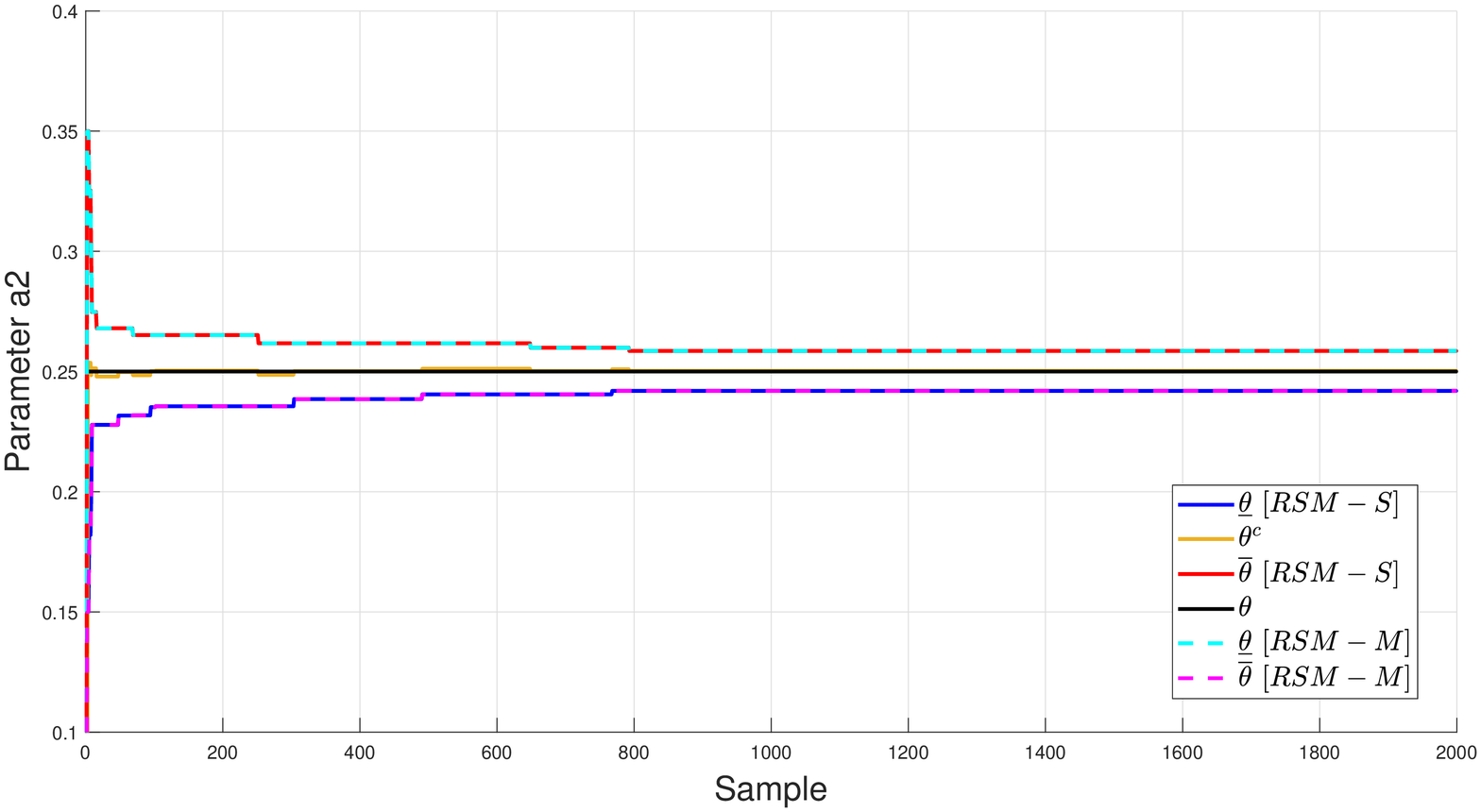}\\
          \includegraphics[width=1\textwidth]{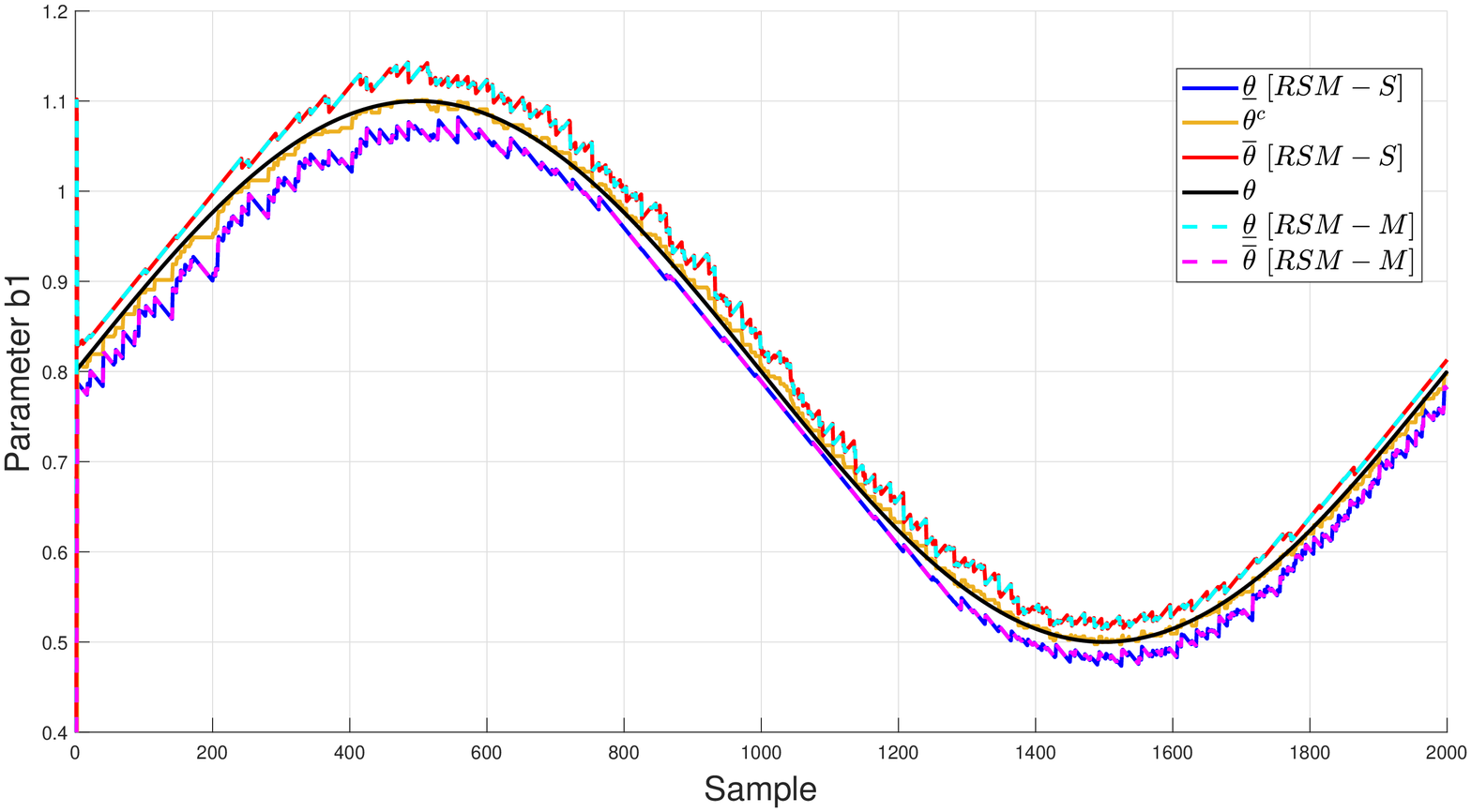}
          \subcaption{$SNR_x = 52$ dB and $SNR_w = 51$.}
          \end{minipage}
           \quad
\begin{minipage}{0.43\textwidth}
\includegraphics[width=1\textwidth]{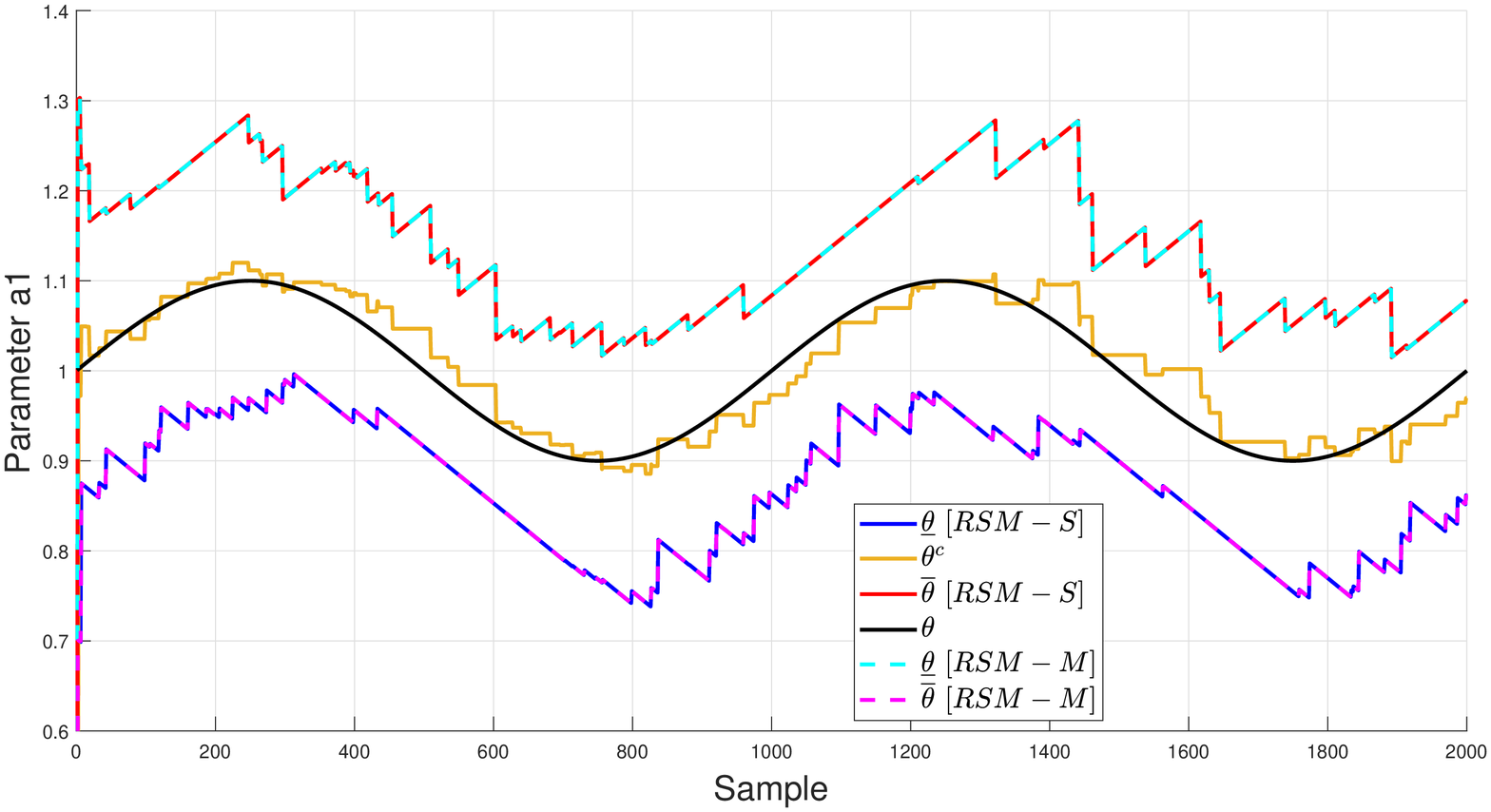}\\
          \includegraphics[width=1\textwidth]{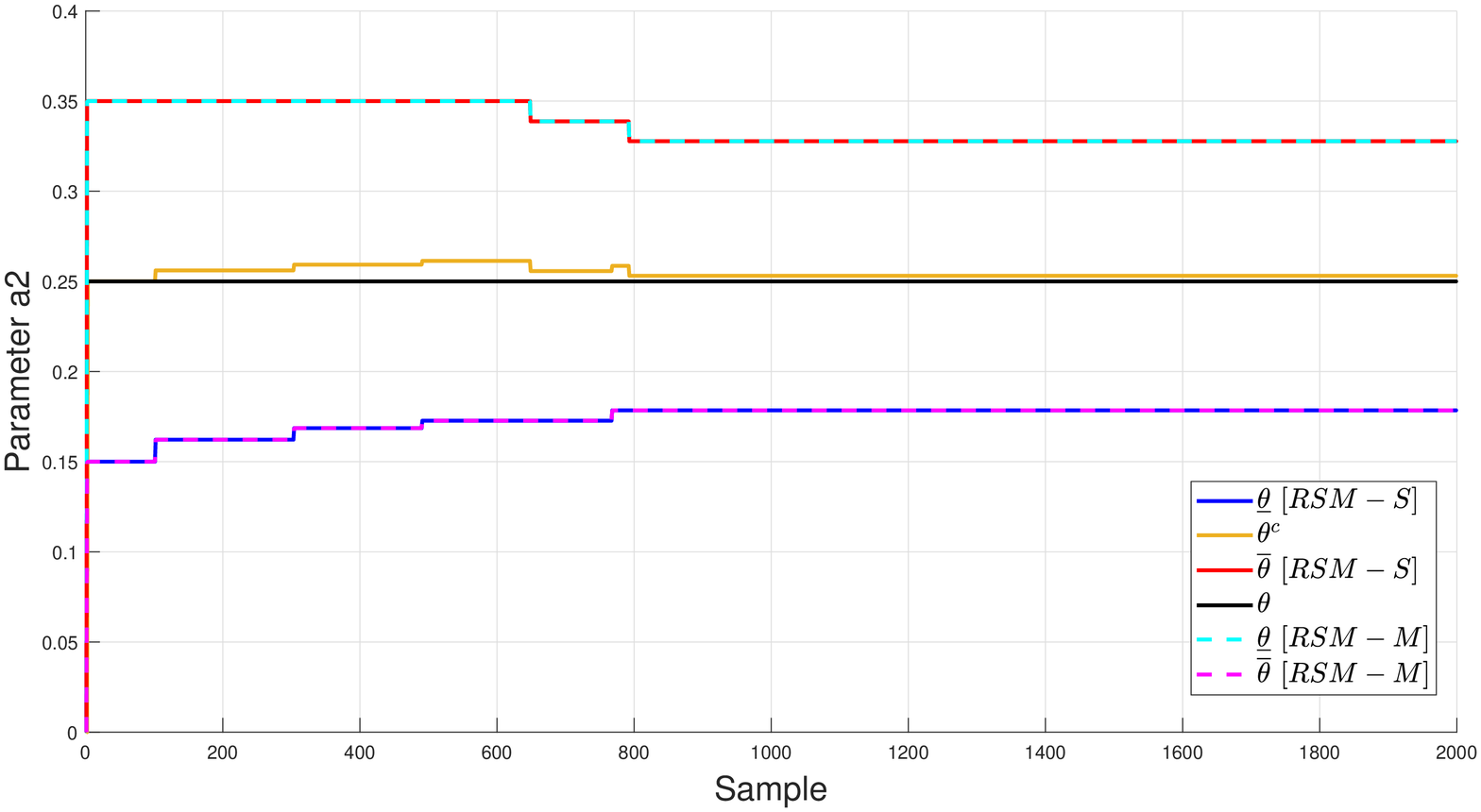}\\
          \includegraphics[width=1\textwidth]{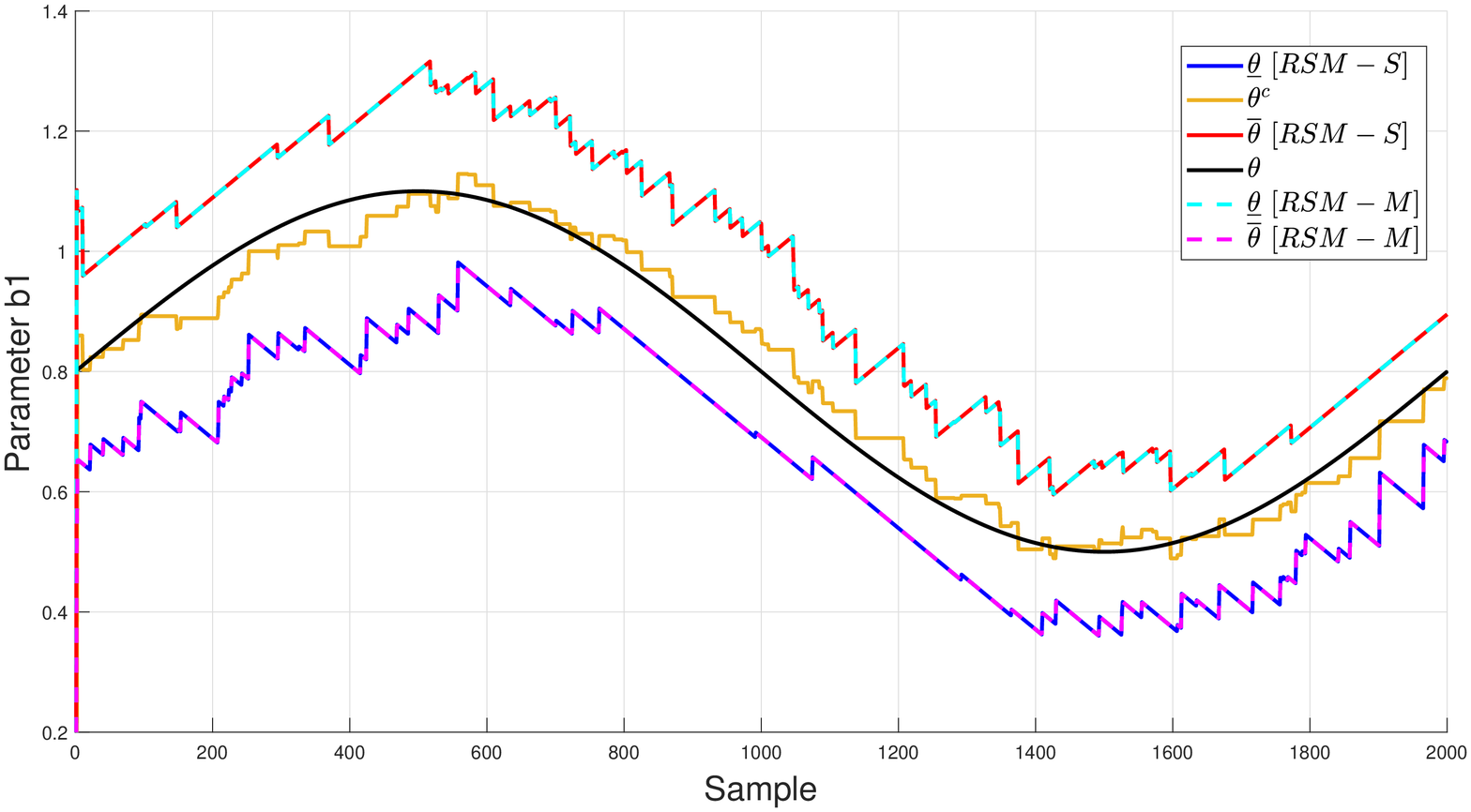}
          \subcaption{$SNR_x = 32$ dB and $SNR_w = 31$.}
          \end{minipage} 
          \caption{Example 2: Comparison between the PUIs computed through $RSM$-$S$ and $RSM$-$M$.}
          \label{fig: 52}
\end{figure*}

The input is a random sequence uniformly distributed between $[-1 \ , \ +1]$. Both input and output sequence are corrupted by random additive noise, uniformly distributed between $[-\Delta_{\zeta} \ , \ \Delta_{\zeta}]$ and $[-\Delta_{\eta} \ , \ \Delta_{\eta}]$, respectively. The following values for the input and output signal-to-noise ratios have been considered in this example: $SNR_x = [52, 32]$ dB and $SNR_w = [51, 31]$ dB. The length of the data set is $N=2000$. In Fig. \ref{fig: 52}, we show a comparison between $\underline{\theta}$, $\overline{\theta}$ and $\theta^{c}$ computed through the algorithm proposed in \cite{cer20}, referred to as recursive set-membership with known signs ($RSM$-$S$), and the one presented in this work, recursive set-membership with McCormick relaxation ($RSM$-$M$). Average elapsed CPU time are quite the same for the two algorithms (1.8 ms for the $RSM$-$S$, and 2.3 ms for $RSM$-$M$).

It can be clearly noticed that the bounds computed through both algorithms are overlapping and are perfectly aligned confirming that the algorithm proposed in this work is able to compute tight PUI (global optimal solution of problem (15)) despite no information on the parameter sign is exploited. 

\clearpage{\pagestyle{empty}\cleardoublepage}
\section{Conclusions}\label{Conclusions}
A novel recursive parameter bounding procedure for SISO discrete-time LTV systems in presence of input and output bounded measurement noise is presented. First the problem is formulated as a nonconvex polynomial optimization problem. Then, based on McCormick envelopes convex relaxation, we show that the parameter uncertainty intervals  for the LTV system can be computed by means of linear programming without assuming any a-priori information on the parameter signs. The effectiveness of the proposed identification scheme is demonstrated by means of two simulation examples.
\bibliographystyle{unsrt}    
\bibliography{Paper_bib}

\end{document}